\documentclass[12pt]{article}
\usepackage{amsmath}
\begin{document} 

\begin{center}
\textsc{Functional Fourier Transformation}
\end{center}

\bigskip
\bigskip

\begin{center}
{\large{\bf INFINITESIMAL FOURIER TRANSFORMATION\\
 FOR THE SPACE OF FUNCTIONALS}}
\end{center}

\bigskip

\begin{center}
\textsc{Takashi Nitta and Tomoko Okada}
\end{center}

\bigskip
\bigskip
                                
\begin{abstract}
The purpose is to formulate a Fourier transformation for the
space of functionals, as an infinitesimal meaning. We extend
${\bf R}$ to
$\,^{\star}(\,^{\ast}{\bf R})$ under the base of nonstandard methods for the
construction. The domain of a functional is the set of all internal functions
from a $\,^{\ast}$-finite lattice to a
$\,^{\ast}$-finite lattice with a double meaning. Considering a $\,^{\ast}$-finite
lattice with a double meaning, we find how to treat the domain for a functional in
our theory of Fourier transformation, and calculate two typical examples.
\end{abstract}

\bigskip
\bigskip
                        
{\bf 0. Introduction}

\medskip

Recently many kinds of geometric invariants are defined on manifolds and they
are used for studying low dimensional manifolds, for example, Donaldson's
invariant, Chern-Simon's invariant and so on. They are originally defined as
Feynman path integrals in physics. The Feynman path integral is in a sense an
integral of a functional on an infinite dimensional space of functions. We would
like to study the Feynman path integral and the originally defined invariants.
For the purpose, we would be sure that it is necessary to construct a theory of
Fourier transformation on the space of functionals. For it, as the later
argument, we would need many stages of infinitesimals and infinites, that is,
we need to put a concept of stage on the field of real numbers. We use
nonstandard methods to develop a theory of Fourier transformation on
the space of functionals.

Feynman([F-H]) used the concept of his path integral for physical quantizations.
The word $''$physical quantizations$''$ has two meanings : one is for quantum
mechanics and the other is for quantum field theory. We usually use the
same word $''$Feynman path integral$''$. However the meanings
included in $''$Feynman path integral$''$ are two sides, according to the above.
One is of quantum mechanics and the other is of quantum field theory. To understand
the Feynman path integral of the first type, Fujiwara([F]) studied it as a
fundamental solution, and Hida([H]), Ichinose, Tamura([Ic],[I-T]) studied it from
their stochastical interests and obtained deep results, using standard
mathematics. In stochastic mathematics, Loeb([Loe]) constructed Loeb measure
theory and investigated Brownian motion that relates to It\^o integral([It]).
Anderson([An]) deloped it. Kamae([Ka]) proved Ergodic theory using nonstandard
analysis. From a nonstandard approach, Nelson([Ne]), Nakamura([Na1],[Na2]) studied
Schr\"odinger equation, Dirac equation and Loo([Loo1],[Loo2]) calculated rigidly
the quantum mechanics of harmonic oscillator. It corresponds to functional
analysis on the space of functions in standard mathematics. 

On the other hand, we would like to construct a frame
of Feynman path integral of the second type, that is, a functional analysis on
the space of functionals. Our idea is the following : in nonstandard analysis,
model theory, especially non-well-founded set theory([N-O-T]), we can extend
${\bf R}$ to $\,^{\ast}{\bf R}$ furthermore a double extension $\,^{\star}(
\,^{\ast}{\bf R})$, and so on. For formulation of Feynman path integral of the
first type, it was necessary only one extension $\,^{\ast}{\bf R}$ of
${\bf R}$ in nonstandard analysis([A-F-HK-L]). In fact there exists an infinite
in
$\,^{\ast}{\bf R}$, however  there are no elements in $\,^{\ast}{\bf
R}$, that is greater than images of the infinite for any functions. The
same situation occurs for infinitesimals. Hence we consider to need a
further extension of
${\bf R}$ to construct a formulation of Feynman path integral of the
second type. If the further extension satisfies some condition, the
extension
$\,^{\star}(\,^{\ast}{\bf R})$ has a higher degree of
infinite and also infinitesimal. We use these to formulate the space of
functionals. We would like to try to construct a theory of Fourier
transformation on the space of functionals and calculate two typical
examples of it.

Historically, for the theories of Fourier transformations in nonstandard
analysis, in 1972, Luxemburg([Lu]) developed a theory of Fourier series with
$\,^{\ast}$-finite summation on the basis of nonstandard analysis. The basic
idea of his approach is to replace the usual ${\infty}$ of the summation to an
infinite natural number $N$. He approximated the Fourier transformation on the
unit circle by the Fourier transformation on the group of
$N$th roots of unity.

Gaishi Takeuti([T]) introduced an infinitesimal delta function $\delta$, and 
Kinoshita([Ki]) defined in 1988 a discrete Fourier transformation for each even
$\,^{\ast}$-finite number $H(\in \,^{\ast}{\bf R}$) :
$(F\varphi)(p)=\sum_{-\frac{H^2}{2}\le z<\frac{H^2}{2}}\frac{1}{H}\exp(-2\pi
ip\frac{1}{H} z)\varphi(\frac{1}{H}z)$, called "infinitesimal Fourier
transformation". He developed a theory for the infinitesimal Fourier
transformation and studied the distribution space deeply, and proved the same
properties hold as usual Fourier transformation of $L^2({\bf R})$. Especially
saying, the delta function $\delta$ satisfies that $\delta^2$, $\delta^2$, ...
, $\sqrt{\delta}$, ... are also hyperfunctions as their meaning, and $F\delta=1$,
$F\delta^2=H$, $F\delta^3=H^2$, ... , $F\sqrt{\delta}=\frac{1}{\sqrt{H}}$, ... .  

In 1989, Gordon([G]) independently defined a generic, discrete
Fourier transformation for each infinitesimal $\Delta$ and
$\,^{\ast}$-finite number $M$, defined
by\\$(F_{\Delta,M}\,\varphi)(p)=\sum_{-M\le z\le M}\Delta\exp(-2\pi ip\Delta 
z)\varphi(\Delta z).$ He studied under which condition the discrete
Fourier transformation
$F_{\Delta,M}$ approximates the usual Fourier transformation
$\mathcal F$ for $L^2({\bf R})$. His proposed condition is ($A'$) of his
notation : let
$\Delta$ be an infinitely small and $M$ an infinitely large natural
number such that $M\cdot\Delta$ is infinitely large. He showed that under the condition
($A'$) the standard part of $F_{\Delta,M}\,\varphi$ approximates the
usual $\mathcal{F}\varphi$ for $\varphi\in L^2({\bf R})$. One of the different
points between Kinoshita's and Gordon's is that there is the term
$\Delta\exp(-2\pi ip\Delta  M)\varphi(\Delta M)$ in the summation of their
two definitions or not. We mention that both definitions are same for the
standard part of the dicrete Fourier transformation for
$\varphi\in L^2({\bf R})$ and Kinoshita's definition satisfies the
condition ($A'$) for an even infinite number
$H$ if $\Delta=\frac{1}{H}$, $M=\frac {H^2}{2}$.

We shall extend their theory of Fourier transformation for the space of
functions to a thery of Fourier transformation for the space of functionals. For
the purpose of this,  we shall represent a space of functions from
${\bf R}$ to ${\bf R}$ as a space of functions from a set of lattices in an
infinite interval $\left[-\frac H2,\frac H2\right)$ to a set of lattices in
an infinite interval $\left[-\frac {H'}{2},\frac {H'}{2}\right)$. We consider
what $H'$ is to treat any function from ${\bf R}$ to ${\bf R}$. If we put a
function $a(x)=x^n(n\in {\bf Z}^+)$, we need that $\frac {H'}{2}$ is greater
than $\left(\frac {H}{2}\right)^n$, and if we  choose a function $a(x)=e^x$, we
need that $\frac {H'}{2}$ is greater than $e^{\frac H2}$. If we choose any
infinite number, there exists a function whose image is beyond the infinite 
number. Since we treat all functions from ${\bf R}$ to
${\bf R}$, we need to put $\frac
{H'}{2}$ as an infinite number greater than any infinite number
of $\,^{\ast}{\bf R}$. Hence we make $\left[-\frac {H'}{2},\frac
{H'}{2}\right)$ not in $\,^{\ast}{\bf R}$ but in
$\,^{\star}(\,^{\ast}{\bf R})$, where $\,^{\star}(\,^{\ast}{\bf R})$ is a
double extension of ${\bf R}$, that is, $H'$ is an infinite number in
$\,^{\star}(\,^{\ast}{\bf R})$. First we shall develop an infinitesimal
Fourier transformation theory for the  space of functionals, and secondly we
calculate foundamental two examples for our infinitesimal Fourier
transformation. In our case, we define an infinitesimal delta function $\delta$
satisfies that $F\delta=1$,
$F\delta^2=H'{}^{H^2}$, $F\delta^3=H'{}^{2H^2}$, ... ,
$F\sqrt{\delta}=H'{}^{-\frac 12 H^2}$, ... , that is, $F\delta^2$,
$F\delta^3$, ... are infinite and $F\sqrt{\delta}$, ... are infinitesimal. These
are a functional
$f$ and an infinite-dimensional Gaussian distribution
$g$ where {\bf st}$(f(\alpha))=\exp\left(\pi
i\int_{-\infty}^{\infty}\alpha^2(t)dt\right),\; ${\bf st}$(g(\alpha))=\exp
\left(-\pi
\int_{-\infty}^{\infty}\alpha^2(t)dt\right)$ for $\alpha\in L^2({\bf R})$. We
obtain the following results of standard meanings : $(Ff)(b)=\overline{f(b)}$ or
$-\overline{f(b)}$ and
$ (Fg)(b)=C_2(b)g(b)$, {\bf st}$(C_2(b))=1$ if $b$ is finite valued. Our
infinitesimal Fourier transformation of $g$ is also $g$ when the domain of $g$ is
standard.

\bigskip 

{\bf 1. Formulation} (cf.[S],[T],[Ki]).

\medskip

To explain our infinitesimal Fourier transformation for the space of
functionals, we introduce Kinoshita's infinitesimal Fourier
transformation for the space of functions. We fix an infinite set
$\Lambda$ and an ultrafilter
$F$ of
$\Lambda$ so that $F$ includes the Fr\'echet filter $F_0(\Lambda)$. We remark
that the set of natural numbers is naturally embedded in $\Lambda$. Let
$H$ be an even infinite number where the definition being even is the following
: if
$H$ is written as
$[(H_{\lambda},\lambda\in\Lambda )]$ then $\{\lambda\in\Lambda\,|\,H_{\lambda}$ 
is even$\}\in F$, where $[\;\;\;]$ denotes the equivalence class with respect to
the ultrafilter $F$. Let
$\varepsilon$ be
$\frac 1H$, that is, if $\varepsilon$ is
$[(\varepsilon_{\lambda},\lambda\in\Lambda)]$ then
$\varepsilon_{\lambda}$ is $\frac{1}{H_{\lambda}}$. Then we shall define
a lattice space ${\bf L}$, a sublattice space
$L$ and a space of functions $R(L)$ :

${\bf L}:=\varepsilon\,^{\ast}{\bf Z}=\{\varepsilon z \,|\, z \in \, ^{\ast}{\bf
Z}\}$,

$L:=\left\{\varepsilon z \,\left|\, z \in \, ^{\ast}{\bf Z},\, -\frac H2 \le
\varepsilon z < \frac H2 \right\}\right.$

$\;\;\;\;\,=\left\{[(\varepsilon_{\lambda}z_{\lambda}),\,\lambda\in\Lambda]\,\left|\,\varepsilon_{\lambda}z_{\lambda}\in
L_{\lambda}\right.\right\}\;\;(\subset {\bf L})$

$R(L):=\left\{\varphi \,\left|\, \varphi \textrm{ is an 
internal
function from } L\textrm{ to }^{\ast}{\bf
C}\right\}\right.$

$\;\;\;\;\;\;\;\;\;\,\,=\left\{[(\varphi_{\lambda},\lambda\in\Lambda)]\,\left|\,\varphi_{\lambda}\textrm{
is a function from }
L_{\lambda}
\textrm{ to }{\bf
C}\right.\right\}$,

\noindent
where $L_{\lambda}:=\left\{
\varepsilon_{\lambda}z_{\lambda}\,\left|\, z_{\lambda} \in
{\bf Z},\, -\frac{H_{\lambda}}{2}\le\varepsilon_{\lambda}z_{\lambda}
<\frac{H_{\lambda}}{2}\right.\right\}$.

Gaishi Takeuti([T]) introduced an infinitesimal delta function $\delta (x)(\in
R(L))$ and Kinoshita([Ki]) defined an infinitesimal Fourier transformation on
$R(L)$. From now on, functions in $R(L)$ are extended to periodic 
functions on ${\bf L}$ with
the period $H$ and we denote them by the same notations. For $\varphi(\in
R(L))$, the infinitesimal Fourier transformation
$F\varphi$, the inverse infinitesimal Fourier transformation $\overline
F\varphi$, and the convolution of
$\varphi$, $\psi(\in R(L))$ are defined as follows : 
\begin{equation*}
\delta (x):= 
\begin{cases}
H & \text{$(x=0)$},\\
0  & \text{$(x \ne 0)$},
\end{cases}
\end{equation*}

$(F\varphi)(p) := \sum _{x \in L} \varepsilon  \exp \left(-2\pi i
px\right)\varphi(x),\;\;(\overline{ F}\varphi)(p) := \sum _{x \in L}
\varepsilon\exp\left(2\pi i px\right)\varphi(x)$,

$(\varphi\ast \psi)(x) :=
\sum _{y \in L} \varepsilon \varphi(x-y)\psi(y)$.

\noindent
He obtained the following	equalities as same as the usual Fourier
analysis :

$\delta = F1=\overline{F} 1,\;
F\textrm{ is unitary},\, F^4 =1, \,\overline{F} F =F\overline{F}= 1$,

$\varphi\ast \delta = \delta \ast \varphi =\varphi,\; 
\varphi \ast \psi = \psi \ast \varphi$,

$F(\varphi \ast \psi )= (F\varphi)(F \psi),\;
F(\varphi \psi) = (F\varphi)\ast (F \psi)$,

$\overline{F}( \varphi \ast  \psi )= (\overline{F}\varphi)(\overline{F}  \psi),\;
\overline{F}(\varphi \psi) = (\overline{F}\varphi)\ast (\overline{F} \psi)$.

\noindent
The most different point is that $\delta^l\;(l\in{\bf R}^{+})$ are also elements
of
$R(L)$ and the Fourier transformation are able to be calculated as
$F\delta^l=H^{(l-1)}$, by the above definition.

On the other hand, we obtain the following theorem from his result 
and an elementary calculation :

\medskip

{\textsc Theorem} 1.1.\;\;{\it For an internal function with two variables
$f:L\times L\to
\,^{\ast}{\bf C}$  and $g(\in R(L))$,

$F_x\left(\sum_{y \in L}\varepsilon f(x-y,y)g(y)\right)(p)
=\left\{F_y(F_u (f(u,y))(p))\ast F_y(g(y))\right\}(p)$,

\noindent
$\textrm{where } F_x,\,F_y,\,F_u \textrm{ are Fourier transformations for }x,\, y,\,
u$, and $\ast$ is the convolution for the variable paired with $y$ by the
Fourier transformation.} 

\medskip

{\it Proof.}\;\; By the above Kinoshita's result, $F(\varphi \psi) =
(F\varphi)\ast (F
\psi)$. We use it and obtain the following :

$F_x\left(\sum_{y \in L}\varepsilon f(x-y,y)g(y)\right)(p)
=\sum_{x,y\in L}\varepsilon \exp(-2\pi ipx) \varepsilon f(x-y,y)g(y)$

$=\sum_{y,u \in L}\varepsilon^2 \exp(-2\pi ip(y+u)) f(u,y)g(y)\;(u:=x-y)$

$=\sum_{y \in L}\left(\varepsilon\exp(-2\pi ipy)
\left(\sum_{u \in L}\varepsilon\exp(-2\pi pu)f(u,y)\right)g(y)\right)$

$=F_y(F_u (f(u,y))(p)\cdot g(y))(p)
=\left\{F_y(F_u (f(u,y))(p))\ast F_y(g(y))\right\}(p)$.

\medskip

To treat a $\,^{\ast}$-unbounded functional $f$ in the nonstandard analysis, we
need a second nonstandardization. Let
$F_2:=F$ be a nonprincipal ultrafilter on an infinite set $\Lambda_2:=\Lambda$ as
above. Denote the ultraproduct of a set $S$ with respect to $F_2$ by $\,^{\ast}S$
as above. Let $F_1$ be another nonprincipal ultrafilter on an infinite set
$\Lambda_1$. Take the $\,^{\ast}$-ultrafilter $\,^{\ast}F_1$ on
$\,^{\ast}\Lambda_1$. For an internal set $S$ in the sense of
$\,^{\ast}$-nonstandardization, let $\,^{\star}S$ be the $\,^{\ast}$-ultraproduct
of $S$ with respect to $\,^{\ast}F_1$. Thus, we define a double ultraproduct
$^{\star}(^{\ast}{\bf R})$, $^{\star}(^{\ast}{\bf Z})$, etc for the set ${\bf
R}$, ${\bf Z}$, etc. It is shown easily that

$$^{\star}(^{\ast}{\bf S})=S^{\Lambda_1\times\Lambda_2}/F_1^{F_2},$$

\medskip
\noindent
where $F_1^{F_2}$ denotes the ultrafilter on $\Lambda_1\times\Lambda_2$ such
that for any $A\subset\Lambda_1\times\Lambda_2$, $A\in F_1^{F_2}$ if and only if 

$$\{\lambda\in\Lambda_1\,|\,\{\mu\in\Lambda_2\,|\,(\lambda,\mu)\in A\}\in
F_2\}\in F_1.$$

\medskip
\noindent
We always work with this double nonstandardization. The natural imbedding
$\,^{\star}S$ of an internal element $S$ which is not considered as a set in
$\,^{\ast}$-nonstandardization is often denoted simply by $S$.

\medskip   

\textsc{Definition} 1.2 (cf.[N-O]).\;\; Let $H (\in \,^{\ast}{\bf Z}),\, H' (\in
\,^{\star}(^{\ast}{\bf Z}))$ be even positive numbers such that $H'$ is larger
than any element in $\,^{\ast}{\bf Z}$, and let $\varepsilon (\in
\,^{\ast}{\bf R}),\,\varepsilon ' (\in \,^{\star}(^{\ast}{\bf R}))$ be 
infinitesimals satifying $\varepsilon H =1,\,\varepsilon ' H' =1$. We define as
follows :

${\bf L} :=\varepsilon \, ^{\ast}{\bf Z}= \{\varepsilon z \,|\,z \in \, 
^{\ast}{\bf
Z}\},\;\;
{\bf L}' :=\varepsilon ' \, ^{\star}(\,^{\ast}{\bf Z})=\{\varepsilon ' z'
\,|\, z' \in \, ^{\star}(\,^{\ast}{\bf Z})\},$
 
$L := \left\{\varepsilon z \,
\left|\,z \in \, ^{\ast}{\bf Z},\, -\frac H2 \le
\varepsilon z < \frac H2 \right\}\right.(\subset {\bf L})$,

$L' :=\left\{\varepsilon' z' \,\left |\,z' \in \, ^{\star}(\,^{\ast}{\bf Z}),\,
-\frac {H'}{2} \le \varepsilon ' z' <
\frac {H'}{2} \right\}(\subset {\bf L}').\right.$

\medskip
\noindent
Here $L$ is an ultraproduct of lattices 

$L_{\mu}:=\left\{
\varepsilon_{\mu}z_{\mu}\,\left|\,z_{\mu}
\in{\bf Z},\, -\frac{H_{\mu}}{2}\le\varepsilon_{\mu}z_{\mu}
<\frac{H_{\mu}}{2}\right.\right\}\;(\mu\in\Lambda_2)$ 

\noindent
in ${\bf R}$, and $L'$ is
also an ultraproduct of lattices 

$L'_{\lambda}
:=\left\{\varepsilon '_{\lambda}z'_{\lambda}\,\left|\, z'_{\lambda}\in \,^{\ast}{\bf
Z},\, -\frac{H'_{\lambda}}{2}\le\varepsilon
'_{\lambda}z'_{\lambda}
<\frac{H'_{\lambda}}{2}\right\}\right.\;(\lambda\in\Lambda_1)$

\noindent
in $\,^{\ast}{\bf R}$ that is an ultraproduct of 

$L'_{\lambda\mu}:=\left\{\varepsilon
'_{\lambda\mu}z'_{\lambda\mu}\,\left|\,z'_{\lambda\mu}
\in{\bf Z},\, -\frac{H'_{\lambda\mu}}{2}\le\varepsilon
'_{\lambda\mu}z'_{\lambda\mu}<\frac{H'_{\lambda\mu}}{2}\right\}
\right.\;(\mu\in \Lambda_2)$. 

\noindent
We define a latticed space of functions $X$ as follows,

$X := \{a \,|\, a \textrm{ is an internal function with double meamings,
from}
\star(L)\textrm{ to } L' \}$

$\;\;\;\;\;=\{[(a_{\lambda}),\,\lambda\in\Lambda_1]\,|\,a_{\lambda}\textrm{ is an
internal function from }L\textrm{ to }
L'_{\lambda}\}$,

\noindent
where $a_{\lambda}:L\to L'_{\lambda}$ is $a_{\lambda}=[(a_{\lambda\mu}),\,\mu\in
\Lambda_2],\; a_{\lambda\mu}: L_{\mu}\to L'_{\lambda\mu}$. 

\noindent
We define three
equivarence relations $\sim_H$, $\sim_{\star(H)}$ and
$\sim_{H'}$ on
${\bf L}$, $\star({\bf L})$ and ${\bf L}'$ :

$x\sim_H y\Longleftrightarrow x-y\in H\,^{\ast}{\bf Z},\;\;
x\sim_{\star(H)} y\Longleftrightarrow x-y\in \star(H)\,^{\star}(\,^{\ast}{\bf Z}),$

$x\sim_{H'} y\Longleftrightarrow x-y\in H'\,^{\star}(\,^{\ast}{\bf Z})$.

\noindent
Then we identify ${\bf L}/\sim_H$, $\star({\bf L})/\sim_{\star(H)}$ and ${\bf
L}'/\sim_{H'}$ as $L$, $\star(L)$ and $L'$. Since $\star(L)$ is identified 
with $L$, the set 
$\star({\bf L})/\sim_{\star(H)}$ is identified with ${\bf L}/\sim_H$. Furthermore
we represent
$X$ as the following internal set : 

\noindent
$\{a\,|\,a\textrm{ is an internal function with double meamings,
from }\star({\bf L})/\sim_{\star(H)}\textrm{to }{\bf
L}'/\sim_{H'}\}.$

\noindent
We use the same notation as a function from $\star(L)$
to $L'$ to represent a function in the above internal set. We define the space $A$
of functionals as follows :

$A:= \{ f \,|\, f \textrm{ is an internal function with double meamings, from
} X 
\textrm{ to } \,^{\star}(\,^{\ast}{\bf C})
\}.$

\noindent
Then $f$ is written as $f=[(f_{\lambda}),\,\lambda\in\Lambda_1]$, $f_{\lambda}$ is
an internal function from the set $\{a_{\lambda}\,|\,a_{\lambda}\textrm{ is an internal
function from } L\textrm{ to }L'_{\lambda}\}$ to $\,^{\ast}{\bf C}$, and
$f_{\lambda}$ is written as $f_{\lambda}=[(f_{\lambda\mu}),\,\mu\in \Lambda_2],\;
f_{\lambda\mu}:
\{a_{\lambda\mu}:L_{\mu}\to L'_{\lambda\mu}\}\to{\bf C}$.

We define an infinitesimal delta function $\delta (a)(\in
A)$, an infinitesimal Fourier transformation of
$f(\in A)$, an inverse infinitesimal Fourier transformation of $f$ and a
convolution of $f$,
$g(\in A)$, by the following :

\medskip

\textsc{Definition} 1.3.
\begin{equation*}
\delta (a):= 
\begin{cases}
(H')^{(\,^{\star}H)^2} & \text{$(a=0)$}, \\
0  & \text{$(a \ne 0)$},
\end{cases}
\end{equation*}

$\varepsilon _0 := (H')^{-{(\,^{\star}H)^2}}
\in \,^{\star}(^{\ast}{\bf R})$,

$(Ff)(b) := \sum _{a \in X} \varepsilon _0  \exp \left(-2\pi i \sum
_{k\in L} a(k)b(k)\right)f(a)$,

$(\overline{ F} f)(b) := \sum _{a \in X}
\varepsilon _0  \exp \left(2\pi i \sum _{k\in L} a(k)b(k)\right)f(a)$,

$(f\ast g)(a)
:=\sum _{a' \in X} \varepsilon _0  f(a-a')g(a')$.

\noindent
We define an inner product on $A$ : $(f,g):=\sum_{b\in
X}\varepsilon_0\overline{f(b)}g(b)$, where $\overline{f(b)}$ is the complex conjugate
of $f(b)$. Then we obtain the following theorem :

\medskip

\textsc{Theorem} 1.4.

$\textrm{(1)}\;\; \delta = F1=\overline{F} 1,\;\;
\textrm{(2)}\;\; F \: \textrm{is unitary}, \, F^4 =1, \overline{ F} F =F
\overline{ F}=
1$,

$\textrm{(3)}\;\; f\ast \delta = \delta \ast f = f ,\;\;
\textrm{(4)}\;\; f \ast g = g \ast f$,

$\textrm{(5)}\;\; F( f \ast g )= (Ff)(Fg),\;\;
\textrm{(6)}\;\; \overline{F}( f \ast g )= (\overline{F} f)(\overline{F} g)$,

$\textrm{(7)}\;\; F(fg) = (Ff)\ast (Fg),\;\;
\textrm{(8)}\;\; \overline{F}(fg) = (\overline{F}f)\ast (\overline{F}g)$.

\medskip

The definition implies the following proposition :

\medskip

\textsc{Proposition} 1.5.\;\; If $l\in{\bf R}^{+}$, then
$F\delta^l=(H')^{(l-1)(\,^{\star}H)^2}$.

\medskip

We define two types of infinitesimal divided differences. Let
$f$ and $a$ be elements of $A$ and $X$ respectively and let
$b (\in X)$ be an internal function whose image is in 
$\,^{\star}(\,^{\ast}{\bf Z})\cap L'$. We remark that $\varepsilon' b$ is an
element of $X$.

\medskip

\textsc{Definition} 1.6.

$(D_{+,b}\, f)(a) := \frac{f(a+ \varepsilon ' b) - f(a)}{\varepsilon '} 
,\;\;
(D_{-,b}\, f)(a) :=\frac{f(a) - f(a- \varepsilon ' b)}{\varepsilon '}$.

\medskip

\noindent
Let $\lambda _{b} (a):=\frac{\exp (2\pi i \varepsilon ' ab)-1}
{\varepsilon '},\;
\overline{\lambda} _{b} (a):= \frac{\exp (-2\pi i \varepsilon '
ab)-1}{\varepsilon '}$. Then we obtain the following theorem corresponding to
Kinoshita's result for the relationship between the infinitesimal Fourier
transformation and the infinitesimal divided differences :

\medskip

\textsc{Theorem} 1.7.

$\textrm{(1)}\;\; (F (D_{+,b} \,f))(a) =\lambda _{b} (a)(Ff)(a),\;\;
\textrm{(2)}\;\; (F (D_{-,b} \,f))(a) = -\overline{\lambda}_{b} (a)(Ff)(a),$

$\textrm{(3)}\;\; (F(\lambda_b f))(a)=-(D_{-,b}\,(Ff))(a),\;\;
\textrm{(4)}\;\;  (F(\overline{\lambda_b} f))(a)= (D_{+,b}\,(Ff))(a),$

$\textrm{(5)}\;\; (D_{+,b} \,(\overline{F}f))(a) =(\overline{F}(\lambda_b f))
(a),\;\;
\textrm{(6)}\;\; (D_{-,b} \,(\overline{F}f))(a) = -(\overline{F}
(\overline{\lambda_b}f))(a),$

$\textrm{(7)}\;\; \lambda _b (a) =2\pi i \left(\frac{\sin (\pi
\varepsilon ' ab)}{\pi \varepsilon '}\right)\exp(\pi i \varepsilon ' ab).$

\medskip

Theorem 1.7 implies the following Corollary :

\medskip
\noindent
\textsc{Corollary} 1.8.\;\; If $\varepsilon' b$ is an
element of $X$, then $(f, D_{+,b}g)=-(D_{+,b}f,g)$ for $f,\, g\in A$.

\medskip

Replacing the definitions of $L'$, $\delta$, $\varepsilon_0$, $F$,
$\overline{F}$ in  Definition 1.2 and Definition 1.3 by the following, we shall
define another type of infinitesimal Fourier transformation. The different point
is only the definition of an inner product of the space of functions $X$. In
Definition 1.3, the
inner product of $a,\,b(\in X)$ is $\sum _{k\in L}a(k)b(k)$, and in the
following definition, it is 
$\,^{\star}\varepsilon\sum _{k\in L}a(k)b(k)$.

\medskip

\textsc{Definition} 1.9.

$L' := \left\{\varepsilon ' z' \,\left|\, z'\in \, ^{\star}(\,^{\ast}{\bf Z}) ,\, -\,
^{\star}H \frac {H'}{2} \le \varepsilon ' z' < \, ^{\star}H \frac {H'}{2}
\right\},\right.$
\begin{equation*}
\delta (a):= 
\begin{cases}
(\,^{\star}H)^{\frac{(\,^{\star}H)^2}{2}}H'{}^{(\,^{\star}H)^2} &
\text{$(a=0)$}, \\ 0  & \text{$(a \ne 0)$},
\end{cases}
\end{equation*}

$\varepsilon_0
:=(\,^{\star}H)^{-\frac{(\,^{\star}H)^2}{2}}H'{}^{-(\,^{\star}H)^2}$

$(Ff)(b) := \sum _{a \in X} \varepsilon _0  \exp
\left(-2\pi i \,^{\star}\!\varepsilon \sum _{k\in L} a(k)b(k)\right)f(a)$,

$(\overline{F} f)(b) :=\sum _{a \in X} \varepsilon _0  \exp \left(2\pi i
\,^{\star}\!\varepsilon \sum _{k\in L} a(k)b(k)\right)f(a).$

\medskip

\noindent
In this case, we obtain the same theorems as Theorem 1.4 and Theorem 1.6, and
the following 
theorem corresponding to Theorem 1.1 :

\medskip

\textsc{Theorem} 1.10.\;\;{\it For an internal function with two variables 
$f:X\times X\to
\,^{\star}(\,^{\ast}{\bf C})$ and $g(\in A)$,

$F_a\left(\sum_{b \in X}\varepsilon_0 f(a-b,b)g(b)\right)(d)=\left\{F_b(F_c
(f(c,b)) (d))\ast F_b (g(b))\right\}(d),$

\noindent
$\textrm{where } F_a,\,F_b,\,F_c \textrm{ are Fourier transformations for }
a,\,
b,\, c$, and $\ast$ is the convolution for the variable pairing with $b$ by the
Fourier transformation.} 

\bigskip

{\bf 2. Proofs of Theorems.}

\medskip

{\it Proof of Theorem} 1.4.

\noindent
(1)
$(F1)(0)=\sum _{a \in X} \varepsilon _0
=\varepsilon _0 ({H'}^2)^{\,(\,^{\star}H)^2}
=(H')^{(\,^{\star}H)^2}$. If $b\ne 0$, then

$(F1)(b)=\sum _{a \in X} \varepsilon _0  \exp 
\left(-2\pi i \sum _{k\in L} a(k)b(k)\right)
=\varepsilon _0\prod_{k \in L}\sum _{a(k) \in L'} \exp(-2\pi ia(k)b(k))$

$=\varepsilon _0\prod_{k \in L,b(k)\ne 0}\sum _{a(k) \in L'} \exp(-2\pi ia(k)b(k))
\cdot\prod_{k \in L,b(k)=0}\sum _{a(k) \in L'} \exp(-2\pi ia(k)b(k))$

$=\prod_{k \in L,b(k)\ne 0}\varepsilon _0 \,\frac{\exp(-2\pi i\varepsilon
'(-\frac{{H'}^2}{2})b(k))(1-\exp(-2\pi i
\varepsilon '{H'}^2 b(k)))}{1-\exp(-2\pi i \varepsilon '
b(k))}$

$\cdot\prod_{k \in L,b(k)= 0}\sum _{a(k) \in L'} \exp(-2\pi ia(k)b(k))
=0.$

\noindent
Hence $F1=\delta$. The same argument implies that $\overline{F}1=\delta$.

\noindent
(2)
$(Ff,Fg)=\sum _{b \in X}\varepsilon _0  \overline{( Ff)(b)}(Fg)(b)$

$=\sum _{b \in X}\varepsilon _0\overline{\sum _{a \in X} \varepsilon _0 
\exp \left(-2\pi i \sum _{k\in L} a(k)b(k)\right)f(a)}\sum _{c \in
X} \varepsilon _0  \exp \left(-2\pi i \sum _{k\in L} c(k)b(k)\right)$

$g(c)$

$=\sum _{a \in X} \sum _{c \in X}\varepsilon _0 ^{2}\overline{
f(a)}g(c)
\sum _{b \in X}\varepsilon _0 \exp \left(-2\pi i \sum _{k\in L} 
(c(k)-a(k))b(k)\right)$

$=\sum _{a \in X} \sum _{c \in X}\varepsilon _0
^{2}\overline{ f(a)}g(c)
\delta (c-a)
=\sum _{a \in X}\varepsilon _0\overline{f(a)}g(a)=(f,g)$.

\noindent
Hence $F$ is unitary. Since $(F^2 f) (c)=(F(Ff))(c)=f (-c)$, $F^4 =1$. Thus
the eigenvalues of $F$ are $1,$ $-1,$ $-i,$ $i$. Furthermore,

$(\overline{F} (Ff))(c)= \sum _{b \in X} \varepsilon _0  \exp \left(2\pi i \sum _{k\in
L} c(k)b(k)\right)\left(\sum _{a \in X} \varepsilon _0  \exp \left(-2\pi i
\sum _{k\in L} a(k)b(k)\right)f(a)\right)$

$=\sum _{a \in X} \left(
\sum _{b \in X} \varepsilon _0^2 
\exp \left(-2\pi i \sum _{k\in L} b(k)(a(k)-c(k))\right)\right)f(a)$

$=\sum _{a \in X} \varepsilon _0\delta (a-c)f(a)=f(c)$.

\noindent
The same argument implies $(F(\overline{F}f))(c)=f(c)$.

\noindent
(3)
$(f\ast \delta)(a)=\sum _{b \in X} \varepsilon _0 f(a-b)
\delta (b)
=f(a)$,

$(\delta \ast f)(a)=\sum _{b \in X} \varepsilon _0 \delta(a-b)
f(b) =f(a)$.

\noindent
(4)
$(f\ast g)(a)=\sum _{b \in X} \varepsilon _0 f(a-b)g(b)
=\sum _{(a-b) \in X} \varepsilon _0 f(a-b)g(a-(a-b))
=(g\ast f)(a).$

\noindent
(5)
$(F(f\ast g))(c) = \sum _{a \in X} \varepsilon _0 \exp \left(-2\pi i \sum
_{k\in L} c(k)a(k) \right) 
\sum _{a \in X} \varepsilon _0 f(a-b)g(b)$

$=\sum _{a \in X} \varepsilon _0 \exp \left(-2\pi i \sum _{k\in L}
c(k)(b(k)+d(k)) \right) 
\sum _{b \in X} \varepsilon _0 f(a-b)g(b)$, where $d(k):=a(k)-b(k)$,

$=\sum _{b \in X} \varepsilon _0 \exp\left(-2\pi i \sum _{k\in L}
c(k)b(k) \right) g(b)
\sum _{d \in X-b} \varepsilon _0 \exp\left(-2\pi i \sum_{k\in L} c(k)d(k)
\right)f(d)$, where $X-b:=\{x-b\,|\,x\in X\}$,

$=\sum _{b \in X} \varepsilon _0 \exp\left(-2\pi i \sum _{k\in L}
c(k)b(k) \right) g(b)
\sum _{d \in X} \varepsilon _0 \exp\left(-2\pi i \sum _{k\in L} c(k)d(k) \right) f(d)
$

$=(Fg)(c)(Ff)(c)
=(Ff)(c)(Fg)(c)$.

\noindent
(6)
Similarly, $\overline F (f' \ast g')=(\overline F f')(\overline F g')$.

\noindent
(7)\;\;The above (6) implies $f' \ast g'=F((\overline F f')(\overline F g'))$. We put
$f'=Ff,\; g'=Fg$. Then we obtain $(Ff) \ast (Fg)=F(fg)$.

\noindent
(8)\;\;Similarly, $(\overline F f) \ast (\overline F g)=\overline F(fg)$.

\medskip

{\it Proof of Theorem} 1.7.

\noindent
(1)$(F(D_{+,b}\,f)) (a) = \sum_{c \in X} \varepsilon _0 \exp (-2\pi i ac)
\frac{1}{\varepsilon '}(f (c+\varepsilon ' b)-f (c))$

$=\sum_{c
\in X} \varepsilon _0 (\frac{1}{\varepsilon '}( \exp (-2\pi i ac) f (c+\varepsilon
' b)-\exp (-2\pi i ac)f (c)))$

$= 
\sum_{c
\in X}
\varepsilon _0 (\frac{1}{\varepsilon '}(\exp (2\pi i \varepsilon ' ab)  (\exp (-2\pi
i a(c+\varepsilon ' b)) f (c+\varepsilon ' b)-\exp (-2\pi i ac) f (c)))$

$=\frac{1}{\varepsilon'}(\exp (2\pi i \varepsilon ' ab) -1)(Ff)(a)=\lambda
_b(a) F f(a)$,

\noindent
(2)
$(F(D_{-,b}\,f)) (a) =\sum_{c \in X} \varepsilon _0 \exp (-2\pi i ac)
\frac{1}{\varepsilon ' }(f (c)- f (c-\varepsilon b ' ))$

$=\sum_{c
\in X} \varepsilon _0 (\frac{1}{\varepsilon ' } 
(\exp (-2\pi i ac) f(c)-\exp (-2\pi i \varepsilon ' ab)
\exp (-2\pi i
a(c-\varepsilon ' b))f (c-\varepsilon ' b)))$

$=\frac{1}{\varepsilon ' }(1-\exp
(-2\pi i \varepsilon ' ab))(F f)(a)=-\overline{\lambda} _b(a) Ff(a)$,

\noindent
(3)
$(F({\lambda} _b f))(a) =\sum_{c \in X} \varepsilon _0  \exp
(-2\pi i ac) ({\lambda} _b f) (c)$

$=\sum_{c \in X} \varepsilon _0 \exp (-2\pi i
ac)\frac{1}{\varepsilon ' }(\exp(2\pi i bc \varepsilon ' )-1)f (c)$

$=\sum_{c \in X} \varepsilon _0 \frac{\exp (-2\pi i (a-b\varepsilon '
)c)-\exp (-2\pi i ac)}{\varepsilon ' } f (c)= -D_{-,b}(Ff) (a)$,

\noindent
(4)
$(F(\overline{\lambda} _b f))(a) =\sum_{c \in X} \varepsilon _0  \exp
(-2\pi i ac) (\overline{\lambda} _b f) (c)$

$=\sum_{c \in X}
\varepsilon _0 \exp (-2\pi i ac)\frac{1}{\varepsilon ' }(\exp(-2\pi i bc
\varepsilon ' )-1)f (c)$

$=\sum_{c \in X} \varepsilon _0 \frac{\exp (-2\pi i (a+b\varepsilon '
)c)-\exp (-2\pi i ac)}{\varepsilon ' } f (c)= D_{+,b}(Ff)(a)$.

\noindent
(1), (2) imply (5), (6).

\medskip

{\it Proof of Theorem} 1.8.

$(f,D_{+,b}\,g)=\sum_{a\in X}\varepsilon _0\overline{f(a)}D_{+,b}\,g
=\sum_{c\in X}\varepsilon _0\overline{(Ff)(c)}(FD_{+,b})g(c)$

$=\sum_{c\in X}\varepsilon _0\overline{(Ff)(c)}\lambda_b(c)(Fg)(c)=\sum_{c\in
X}\varepsilon _0\overline{\overline{\lambda_b(c)}(Ff)(c)}(Fg)(c)$

$=-\sum_{c\in X}\varepsilon _0\overline{F(D_{-,b}\,f)(c)}(Fg)(c)=-\sum_{a\in
X}\varepsilon _0\overline{D_{-,b}\,f(a)}g(a)=-(D_{-,b}\,f,g)$.

\bigskip

{\bf 3. Examples.}

\medskip

We calculate two
examples of the infinitesimal Fourier transformation for the space $A$ of
functionals. Let $\star\circ\ast : {\bf R}\to \,^{\star}(\,^{\ast}{\bf R})$ be the natural
elementary embedding and let ${\bf st}(c)$ for $c\in\,^{\star}(\,^{\ast}{\bf
R})$ be the standard part of $c$ with respect to the natural elementary
embedding $\star\circ\ast$. The first is for
$\exp(i\pi\,^{\star}\varepsilon\sum_{k\in L} a^2(k))$ and the second is for
$\exp(-\pi\,^{\star}\varepsilon\sum_{k\in L} a^2(k))$. We denote the two
functionals by $f(a), g(a)$. If there is an
$L^2$-function $\alpha(t)$ on ${\bf R}$ for $a(k)$ so that
$a(k)=\star((\,^{\ast}\alpha)(k))$, then {\bf st}$(f(a))=\exp\left(i\pi
\int_{-\infty}^{\infty}\alpha^2(t)dt\right)$, and {\bf st}$(g(a))=\exp
\left(-\pi
\int_{-\infty}^{\infty}\alpha^2(t)dt\right)$. Then we obtain the following
results :

\medskip

{\it Example} 1. $(Ff)(b)=C_1\overline{f(b)}$, where $C_1=\sum_{a\in
X}\varepsilon_0\exp(i\pi\,^{\star}\varepsilon\sum_{k\in L}a^2(k))$, it is just
a standard number $(-1)^{\frac H2}$.

\medskip

{\it Example} 2. $(Fg)(b)=C_2(b) g(b)$, where $C_2(b)=\sum_{a\in
X}\varepsilon_0\exp(-\pi\,^{\star}\varepsilon\sum_{k\in L}(a(k)+ib(k))^2)$, and
if $b$ is a finite valued function then it satisfies
that $\textrm{st}\left(\textrm{st}\left(C_2(b)\right)
\right)=1.$

\medskip

For it, we calculate
Kinoshita's infinitesimal Fourier transformation $\varphi_1(x)=\exp(i\pi x^2)$,
$\varphi_2(x)=\exp(-\pi x^2)$ for the space $R(L)$ of functions. We
obtain : 

$(F\varphi_1)(p)=\exp\left(i\frac{\pi}{4}\right)\overline{\varphi_1(p)}$,

$(F\varphi_2)(p)=c(p){\varphi_2(p)}$, where
$\textrm{st}\left(c(p)\right)=\int_{-\infty}^{\infty}\exp(-\pi t^2)dt$, in the
case of finite $p$.

We denote the following :
 
$R({\bf L}):=\{\varphi'\,|\,\varphi'\textrm{ is an internal function from }{\bf
L}\textrm{ to }\,^{\ast}{\bf C}\}$,

$R_H({\bf
L}):=\{\varphi'\in R({\bf L})\,|\,\varphi'(x+H)=\varphi'(x)\}$.

\noindent
Let $e$ be a mapping from $R(L)$ to $R_H({\bf L})$ defined by
$(e(\varphi))(x)=\varphi(\hat x)$, where $\hat x$ is an element of $L$ satisfying $x
\sim_H \hat x$. Now $\exp(i\pi x^2)$ is an element of $R_H({\bf L})$, in
fact, putting $x\in {\bf L}$ :

$\exp(i\pi (x+H)^2)=\exp(i\pi(x^2 +2xH+H^2))$

$=\exp(i\pi x^2)\exp(2\pi i\varepsilon zH)\exp(i\pi H^2)\;\;(x=\varepsilon z (z\in
\,^{\ast}{\bf Z}))$

$=\exp(i\pi x^2)\;\textrm{, as }\varepsilon H=1\textrm{ and }H
\textrm{ is even}$.

\noindent
Hence $e(\exp(i\pi x^2))=\exp(i\pi x^2)$, that is,
$e(\varphi_1(x))=\varphi_1(x)$. We do the infinitesimal Fourier transformation
of $\varphi_1(x)$,

$(F\varphi_1)(p)=
\sum_{x\in L}\varepsilon\exp(i\pi x^2)\exp(-2\pi ixp)
=\sum_{x\in L}\varepsilon\exp(i\pi (x-p)^2)\exp(-i\pi p^2)$

$=\sum_{x-y\in L}\varepsilon\exp(i\pi (x-p)^2)\exp(-i\pi p^2)
=\sum_{x\in L}\varepsilon\exp(i\pi x^2)\exp(-i\pi p^2)$

$=\left(\sum_{x\in L}\varepsilon\exp(i\pi x^2)\right)\overline{\exp(i\pi p^2)}
=\left(\sum_{x\in L}\varepsilon\exp(i\pi
x^2)\right)\overline{\varphi_1(p)}.$

\noindent
By Gauss sums (cf.[R], p.409) : $\sum_{n=0}^{N-1}\exp\left(\frac{2\pi
i}{N}n^2\right)=\frac{1+(-i)^N}{1-i}\sqrt{N},$ when
$N=4m\;(m\in{\bf N})$,
$\sum_{n=0}^{N-1}\exp\left(\frac{2\pi i}{N}n^2\right)=(1+i)\sqrt{N}$. Using it,
we obtain the following :

$\sum_{x\in L}\varepsilon\exp(i\pi
x^2)=\sum_{z=-\frac{H^2}{2}}^{\frac{H^2}{2}-1}\varepsilon\exp(i\pi(\varepsilon
z)^2)=\sum_{z=-\frac{H^2}{2}}^{\frac{H^2}{2}-1}\varepsilon\exp\left(2\pi i
\frac{z^2}{2H^2}\right)$

$=\frac{1}{2}\varepsilon\sum_{z=0}^{2H^2-1}\exp\left(2\pi i
\frac{z^2}{2H^2}\right)=\frac{1}{2}\varepsilon
(1+i)\sqrt{2H^2}=\frac{1+i}{\sqrt{2}}=\exp\left(i\frac{\pi}{4}\right)$.

\noindent
Nextly we calculate the infinitesimal Fourier transformation of
$e(\varphi_2(x))$.

$(F(e(\varphi_2)))(p)=\sum_{x\in L}\varepsilon e(\exp(-\pi x^2))\exp(-2\pi ixp)$

$=\left(\sum_{x\in L}\varepsilon\exp(-\pi (x+ip)^2)\right){\exp(-\pi p^2)}
=\left(\sum_{x\in L}\varepsilon\exp(-\pi (x+ip)^2)\right){\varphi_2(p)}$.

\noindent
We assume that $p(\in L)$ is finite. Since $\exp(-\pi (x+ip)^2)$ is proved to
be an
$S$-integrable function directly, the term
$\sum_{x\in L}\varepsilon\exp(-\pi (x+ip)^2)$ satisfies the following
(cf.[An], [Loe]) :

$\textrm{st}\left(\sum_{x\in L}\varepsilon
\exp(-\pi (x+ip)^2)\right)=\int_{-\infty}^{\infty}\exp(-\pi (t+i \,^{\circ}p)^2)dt$,

\noindent
where $p\in L$, $\,^{\circ}p=\textrm{st}(p)(\in {\bf R})$. We remark that the
integral $\int_{-\infty}^{\infty}\exp(-\pi (t+i \,^{\circ}p)^2)dt$ is equal to
$\int_{-\infty}^{\infty}\exp(-\pi t^2)dt$.

We define an equivarent relation $\sim_{\,^{\star}HH'}$ in ${\bf
L}'$ by $x\sim_{\,^{\star}HH'} y\Leftrightarrow x-y\in
\,^{\star}HH'\,^{\star}(\,^{\ast}{\bf Z})$. We identify ${\bf
L}'/\sim_{\,^{\star}HH'}$ with
$L'$. Let

$X_{H,\,^{\star}HH'}:=\{a'\,|\,a'\textrm{ is an internal function with double
meamings, from }\star({\bf L})/\sim_{\star(H)}\textrm{ to }\; {\bf
L}'/\sim_{\,^{\star}HH'}\}$, 

\noindent
and let ${\bf e}$ is a mapping from $X$ to
$X_{H,\,^{\star}HH'}$, defined by $({\bf e}(a))([k])=[a(\hat k)]$, where
$[\;\;\;]$ in left hand side represents the equivarent class for the
equivarent relation $\sim_{\star(H)}$ in $\star({\bf L})$, $\hat k$ is a
representative in
$\star(L)$ satisfying $k \sim_{\star(H)} \hat k$, and $[\;\;\;]$ in right hand
side represents the equivarent class for the equivarent relation
$\sim_{\,^{\star}HH'}$ in ${\bf L}'$. We consider an example $f(a)=\exp
\left(i\pi \,^{\star}\varepsilon\sum_{k\in L} a^2(k)\right)$ in the space $A$ of
functionals, for
$a\in X$. 

$\exp \left(i\pi
\,^{\star}\varepsilon\sum_{k\in L} (a(k)+\,^{\star}HH')^2\right)$

$=\exp \left(i\pi\,^{\star}\varepsilon\sum_{k\in
L} a^2(k)\right) \exp \left(2i\pi\,^{\star}\varepsilon\sum_{k\in
L} \,^{\star}HH'a(k)\right)\exp \left(i\pi\,^{\star}\varepsilon\sum_{k\in
L} \,^{\star}H^2H'{}^2\right)$

$=\exp \left(i\pi\,^{\star}\varepsilon\sum_{k\in
L} a^2(k)\right)$.

\noindent
Hence if ${\bf e}^{\sharp}(f)(a)$ is defined by $f({\bf e}(a))$, then
${\bf e}^{\sharp}(f)=f$. We do the infinitesimal Fourier transformation of
$f(a)$.

$(Ff)(b)=(F\left(\exp \left(i\pi\,^{\star}\varepsilon\sum_{k\in
L}a^2(k)\right)\right))(b)$

$=\sum_{a\in X}\varepsilon_0\exp \left(-2i\pi
\,^{\star}\varepsilon\sum_{k\in L}a(k)b(k)\right)\exp
\left(i\pi\,^{\star}\varepsilon\sum_{k\in L}a^2(k)\right)$

$=\left(\sum_{a\in X}\varepsilon_0\exp
\left(i\pi\,^{\star}\varepsilon
\sum_{k\in L}(a(k)-b(k))^2\right)\right)\exp \left(-i\pi\,^{\star}\varepsilon 
\sum_{k\in L}b^2(k)\right)$.

\noindent
The term $\sum_{a\in X}\varepsilon_0\exp \left(i\pi\,^{\star}\varepsilon
\sum_{k\in L}(a(k)-b(k))^2\right)$ is represented as 

$\sum_{a_{\lambda\mu}\in X_{\lambda\mu}}(\varepsilon_0)_{\lambda\mu}\exp
(i\pi\varepsilon_{\mu}\sum_{k_{\mu}\in
L_{\mu}}((a(k))_{\lambda\mu}-(b(k))_{\lambda\mu})^2)$

$=\prod_{k_{\mu}\in
L_{\mu}}(\sum_{((a(k))_{\lambda\mu}
\in L_{\lambda\mu}'}\sqrt{\varepsilon_{\mu}}\varepsilon_{\lambda\mu}'\exp
(i\pi\varepsilon_{\mu} ((a(k))_{\lambda\mu}-(b(k))_{\lambda\mu})^2))$

\noindent
for $\lambda\mu$ component, as
$\varepsilon_0=(\,^{\star}\varepsilon)^{\frac{(\,^{\star}H)^2}{2}}
\varepsilon'{}^{(\,^{\star}H)^2}$.
In the above, 
$X_{\lambda\mu}$ is the set of
functions from $L_{\mu}$ to $L_{\lambda\mu}'$. If we put 
$(a(k))_{\lambda\mu}=\varepsilon_{\lambda\mu}'z_{\lambda\mu}^a\;
(z_{\lambda\mu}^a\in
{\bf Z})$ and $(b(k))_{\lambda\mu}=\varepsilon_{\lambda\mu}'z_{\lambda\mu}^b\;
(z_{\lambda\mu}^b\in {\bf Z})$, we remark that $z_{\lambda\mu}^a$ and
$z_{\lambda\mu}^b$ depend on $k_{\mu}$. The above is 

$\prod_{k_{\mu}\in L_{\mu}}(\sum_{\varepsilon_{\lambda\mu}'z_{\lambda\mu}^a \in
L_{\lambda\mu}'}\sqrt{\varepsilon_{\mu}}\varepsilon_{\lambda\mu}'\exp
(i\pi\varepsilon_{\mu}
(\varepsilon_{\lambda\mu}'z_{\lambda\mu}^a-\varepsilon_{\lambda\mu}'
z_{\lambda\mu}^b)^2))$

$=\prod_{k_{\mu}\in L_{\mu}}
(\sum_{\varepsilon_{\lambda\mu}'z_{\lambda\mu}^a \in
L_{\lambda\mu}'}\sqrt{\varepsilon_{\mu}}\varepsilon_{\lambda\mu}'\exp (i\pi
(\sqrt{\varepsilon_{\mu}}\varepsilon_{\lambda\mu}'z_{\lambda\mu}^a
-\sqrt{\varepsilon_{\mu}}\varepsilon_{\lambda\mu}'z_{\lambda\mu}^b)^2))$.

\noindent
Since $\exp(i\pi\varepsilon
x^2)$ is a periodic function with period $\,^{\star}HH'$ on ${\bf L}'$, 
we obtain 

$\sum_{\varepsilon_{\lambda\mu}'z_{\lambda\mu}^a \in
L_{\lambda\mu}'}\sqrt{\varepsilon_{\mu}}\varepsilon_{\lambda\mu}'\exp \left(i\pi
(\sqrt{\varepsilon_{\mu}}\varepsilon_{\lambda\mu}'z_{\lambda\mu}^a
-\sqrt{\varepsilon_{\mu}}\varepsilon_{\lambda\mu}'z_{\lambda\mu}^b)^2\right)$

$=\sum_{\varepsilon_{\lambda\mu}'z_{\lambda\mu}^a \in L_{\lambda\mu}'}
\sqrt{\varepsilon_{\mu}}\varepsilon_{\lambda\mu}'\exp \left(i\pi
(\sqrt{\varepsilon_{\mu}}\varepsilon_{\lambda\mu}'z_{\lambda\mu}^a)^2\right)$.

\noindent
Hence 

\noindent
$\sum_{a\in X}\varepsilon_0\exp \left(i\pi\,^{\star}\varepsilon
\sum_{k\in L}(a(k)-b(k))^2\right)
=\sum_{a\in X}\varepsilon_0\exp \left(i\pi\,^{\star}\varepsilon
\sum_{k\in L}(a(k))^2\right)=C_1$.

We calculate $C_1$ as follows :

$\sum_{a\in X}\varepsilon_0\exp \left(i\pi\,^{\star}\varepsilon
\sum_{k\in L}(a(k))^2\right)
=\varepsilon_0\prod_{x\in L}\sum_{a(k)\in L'}\exp \left(i\pi\,^{\star}
\varepsilon(a(k))^2\right)$

$=\prod_{x\in
L}(\sqrt{\,^{\star}H}H')^{-1}\sum_{z'=-\frac{\,^{\star}HH'{}^{2}}{2}}^{
\frac{\,^{\star}HH'{}^{2}}{2}-1}\exp
\left(i\pi\,^{\star}\varepsilon (\varepsilon'z')^2\right)$

$=\prod_{x\in
L}(\sqrt{\,^{\star}H}H')^{-1}\sum_{z'=-\frac{\,^{\star}HH'{}^{2}}{2}}^{
\frac{\,^{\star}HH'{}^{2}}{2}-1}\exp
\left(i\pi\frac{1}{\,^{\star}H} \frac{1}{H'{}^2}z'{}^2\right)$
 
$=\prod_{x\in
L}(\sqrt{\,^{\star}H}H')^{-1}\frac{1}{2}\sum_{z'=0}^{2\,^{\star}HH'{}^{2}-1}\exp
\left(2\pi i\frac{z^2}{2\,^{\star}HH'{}^2}\right)$

$\left(\textrm{Gauss sums (cf. [R], p.409)}\;
\sum_{n=0}^{N-1}\exp\left(\frac{2\pi
i}{N}n^2\right)=(1+i)\sqrt{N}\;\;(N=4m\;(m\in{\bf N})\right)$

$=\prod_{x\in
L}(\sqrt{\,^{\star}H}H')^{-1}\frac{1}{2}(1+i)\sqrt{2\,^{\star}HH'{}^2}
=\prod_{x\in L}\frac{1+i}{\sqrt{2}}
=\left(\frac{1+i}{\sqrt{2}}\right)^{H^2}
=\left(\exp\left(i\frac{\pi}{4}\right)\right)^{H^2}$

$=\left(\exp\left(i\pi\right)\right)^{\left(\frac{H}{2}\right)^2}
=(-1)^{\left(\frac{H}{2}\right)^2} =(-1)^{\frac{H}{2}}$.

We do the infinitesimal Fourier
transformation of $g(a)$.

$(Fg)(b)=(F\left(\exp \left(-\pi\,^{\star}\varepsilon\sum_{k\in
L}a^2(k)\right)\right))(b)$

$=\sum_{a\in X}\varepsilon_0\exp \left(-2i\pi
\,^{\star}\varepsilon\sum_{k\in L}a(k)b(k)\right)\exp
\left(-\pi\,^{\star}\varepsilon\sum_{k\in L}a^2(k)\right)$

$=\left(\sum_{a\in X}\varepsilon_0\exp \left(-\pi\,^{\star}
\varepsilon\sum_{k\in L}(a(k)+i b(k))^2\right)\right)\exp
\left(-\pi \,^{\star}\varepsilon \sum_{k\in L}b^2(k)\right)$.

\noindent
We consider the term $\sum_{a\in X}\varepsilon_0\exp
\left(-\,^{\star}\varepsilon\pi
\sum_{k\in L}(a(k)+ i b(k))^2\right).$ We write

$(a(k))_{\lambda\mu}=\varepsilon_{\lambda\mu}'z_{\lambda\mu}^a
(k_{\mu})\;\;(z_{\lambda\mu}^a (k_{\mu})\in {\bf
Z})$, $(b(k))_{\lambda\mu}=\varepsilon_{\lambda\mu}'z_{\lambda\mu}^b
(k_{\mu})\;\;(z_{\lambda\mu}^b (k_{\mu})\in {\bf Z}).$

\noindent
From now on, we denote $z_{\lambda\mu}^a (k_{\mu})$, $z_{\lambda\mu}^b 
(k_{\mu})$ by $z_{\lambda\mu}^a$, $z_{\lambda\mu}^b$ for simplicity. Then the
$\lambda\mu$-component of $\sum_{a\in X}\varepsilon_0\exp
\left(-\pi\,^{\star}\varepsilon\sum_{k\in L}(a(k)+ i b(k))^2\right)$ is equal
to 

$\sum_{a_{\lambda\mu}\in X_{\lambda\mu}}(\varepsilon_0)_{\lambda\mu} \exp
\left(-\pi\varepsilon_{\mu}\sum_{k_{\mu}\in L_{\mu}}((a(k))_{\lambda\mu}+ i
(b(k))_{\lambda\mu})^2\right)$

$=\prod_{k_{\mu}\in
L_{\mu}}\left(\sum_{(a(k))_{\lambda\mu}\in
L_{\lambda\mu}'}\sqrt{\varepsilon_{\mu}}\varepsilon'{}_{\lambda\mu} \exp
\left(-\pi\varepsilon_{\mu} ((a(k))_{\lambda\mu}+ i
(b(k))_{\lambda\mu})^2\right)\right)$

$=\prod_{k_{\mu}\in
L_{\mu}}\left(\sum_{\varepsilon_{\lambda\mu}'z_{\lambda\mu}^a \in
L_{\lambda\mu}'}\sqrt{\varepsilon_{\mu}}\varepsilon'{}_{\lambda\mu} \exp
\left(-\pi\varepsilon_{\mu} (\varepsilon_{\lambda\mu}'z_{\lambda\mu}^a+ i
\varepsilon_{\lambda\mu}'z_{\lambda\mu}^b)^2\right)\right)$

$=\prod_{k_{\mu}\in
L_{\mu}}\left(\sum_{\varepsilon_{\lambda\mu}'z_{\lambda\mu}^a
\in L_{\lambda\mu}'}\sqrt{\varepsilon_{\mu}}\varepsilon'{}_{\lambda\mu} \exp
\left(-\pi (
\sqrt{\varepsilon_{\mu}}\varepsilon_{\lambda\mu}'z_{\lambda\mu}^a+ i 
\sqrt{\varepsilon_{\mu}}\varepsilon_{\lambda\mu}'z_{\lambda\mu}^b)
^2\right)\right).$

\noindent
We assume that $b\;(\in X)$ is finitely valued, that is, $\exists b_0\in {\bf
R}\textrm{ s.t. }k\in L\Rightarrow |b(k)|\le
\star(\ast(b_0)).$ The $\lambda\mu$-component of $\displaystyle\frac{\sum_{a\in
X}\varepsilon_0\exp
\left(-\pi\,^{\star}\varepsilon\sum_{k\in L}(a(k)+ i
b(k))^2\right)}{\star\left(\ast\left(\int_{-\infty}^{\infty}\exp(-\pi
x^2)dx\right)^{H^2}\right)}$ is equal to

$\displaystyle\prod_{k_{\mu}\in
L_{\mu}}\frac{\sum_{\varepsilon_{\lambda\mu}'z_{\lambda\mu}^a
\in L_{\lambda\mu}'}\sqrt{\varepsilon_{\mu}}\varepsilon'{}_{\lambda\mu} \exp
\left(-\pi (
\sqrt{\varepsilon_{\mu}}\varepsilon_{\lambda\mu}'z_{\lambda\mu}^a+ i 
\sqrt{\varepsilon_{\mu}}\varepsilon_{\lambda\mu}'z_{\lambda\mu}^b)^2\right)}
{\int_{-\infty}^{\infty}\exp(-\pi
x^2)dx}.$

\noindent
We write
$B_{\lambda\mu}(k_{\mu}):=\sum_{\varepsilon_{\lambda\mu}'z_{\lambda\mu}^a\in
L_{\lambda\mu}'}\sqrt{\varepsilon_{\mu}}\varepsilon'{}_{\lambda\mu} \exp
\left(-\pi (
\sqrt{\varepsilon_{\mu}}\varepsilon_{\lambda\mu}'z_{\lambda\mu}^a +i
\sqrt{\varepsilon_{\mu}}\varepsilon_{\lambda\mu}'z_{\lambda\mu}^b)^2\right)$

$-\int_{-\infty}^{\infty}\exp(-\pi (x+i\sqrt{\varepsilon_{\mu}}
\varepsilon_{\lambda\mu}'z_{\lambda\mu}^b)^2)dx.$

\noindent
It is equal to

$-2\int_{-\infty}^{-\sqrt{\varepsilon_{\mu}}\frac{H_{\lambda\mu}'}{2}}\exp(-\pi
(x+ib_{\lambda\mu})^2)dx$

$+\sum_{\varepsilon_{\lambda\mu}'z_{\lambda\mu}^a \in
L_{\lambda\mu}'}\sqrt{\varepsilon_{\mu}}\varepsilon'{}_{\lambda\mu} \exp
\left(-\pi (
\sqrt{\varepsilon_{\mu}}\varepsilon_{\lambda\mu}'
z_{\lambda\mu}^a+i\sqrt{\varepsilon_{\mu}}\varepsilon_{\lambda\mu}'
z_{\lambda\mu}^b)^2\right)$

$-\int_{-\sqrt{\varepsilon_{\mu}}\frac{H_{\lambda\mu}'}{2}}^{\sqrt{\varepsilon_{\mu}}\frac{H_{\lambda\mu}'}{2}}\exp(-\pi
(x+ib_{\lambda\mu})^2)dx.\;\;\;\;\;\cdots (\ast_1)$

\noindent
Then the above is equal to

$\displaystyle\prod_{k_{\mu}\in
L_{\mu}}\frac{\sum_{\varepsilon_{\lambda\mu}'z_{\lambda\mu}^a
\in L_{\lambda\mu}'}\sqrt{\varepsilon_{\mu}}\varepsilon'{}_{\lambda\mu} \exp
\left(-\pi (
\sqrt{\varepsilon_{\mu}}\varepsilon_{\lambda\mu}'z_{\lambda\mu}^a+ i 
\sqrt{\varepsilon_{\mu}}\varepsilon_{\lambda\mu}'z_{\lambda\mu}^b)^2\right)}
{\int_{-\infty}^{\infty}\exp(-\pi
x^2)dx}$

$=\displaystyle\prod_{k_{\mu}\in
L_{\mu}}\left(1+\frac{B_{\lambda\mu}(k_{\mu})}{\int_{-\infty}^{\infty}\exp(-\pi
x^2)dx}\right)=\displaystyle\left(1+\frac{B_{\lambda\mu}(k_{\mu})}{\int_{-\infty}^{\infty}\exp(-\pi
x^2)dx}\right)^{H_{\mu}^2}$

$=\displaystyle\left(\left(1+\frac{B_{\lambda\mu}(k_{\mu})}{\int_{-\infty}^{\infty}\exp(-\pi
x^2)dx}\right)^{\frac{1}{B_{\lambda\mu}(k_{\mu})}}\right)
^{B_{\lambda\mu}(k_{\mu})H_{\mu}^2}$

$=\displaystyle\left(\left(1+\frac{1}{\int_{-\infty}^{\infty}\exp(-\pi
x^2)dx\cdot\frac{1}{B_{\lambda\mu}(k_{\mu})}}\right)
^{\frac{1}{B_{\lambda\mu}(k_{\mu})}}\right)^{B_{\lambda\mu}(k_{\mu})H_{\mu}^2}.
\;\;\;\;\;\cdots (\ast_2)$

We show that $[(B_{\lambda\mu}(k_{\mu}))]$ is infinitesimal in 
$\,^{\star}(\,^{\ast}{\bf C})$ with respect to ${\bf C}$. It implies that
$\left[\left(\frac{1}{B_{\lambda\mu}(k_{\mu})}\right)\right]$ is infinite in
$\,^{\star}(\,^{\ast}{\bf C})$. For a sequence $a_n$, we
remark that
$\lim_{n\to \infty}a_n=a\Longleftrightarrow \forall 
N: \textrm{infinite with respect to }{\bf C}\;(\in \,^{\star}(\,^{\ast}{\bf
C})),\; \textrm{{\bf st}}(\,^{\star}(\,^{\ast}a_N))=a.$
Hence 

$\displaystyle{\bf
st}\left(\left[\left(1+\frac{1}{\int_{-\infty}^{\infty}\exp(-\pi
x^2)dx\cdot\frac{1}{B_{\lambda\mu}(k_{\mu})}}\right)^{\frac{1}{B_{\lambda\mu}
(k_{\mu})}}\right]\right)=\exp\left(-\int_{-\infty}^{\infty}\exp(-\pi
x^2)dx\right)$

\noindent
and ${\bf st}([(\ast_1))])=1.$ 

Since
$b_{\lambda\mu}$ is finite and
$\left[\left(\sqrt{\varepsilon_{\mu}}\frac{H_{\lambda\mu}'}{2}
\right)\right]$ is infinitesimal in $\,^{\star}(\,^{\ast}{\bf
R})$ with respect to ${\bf R}$, the first term of ($\ast_1$) is infinitesimal 
 in $\,^{\star}(\,^{\ast}{\bf
C})$ with respect to ${\bf C}$. In order to show that
$[(B_{\lambda\mu}(k_{\mu}))]$  is infinitesimal in $\,^{\star}(\,^{\ast}{\bf
C})$, we consider the second and third terms in $(\ast_1)$, and we prove that it
is represents an infinitesimal number. First we calculate 

$\exp(-\pi
(x+i\sqrt{\varepsilon_{\mu}}\varepsilon_{\lambda\mu}'z_{\lambda\mu}^b)^2)-\exp(-\pi
(\sqrt{\varepsilon_{\mu}}\varepsilon_{\lambda\mu}'
z_{\lambda\mu}^a+i\sqrt{\varepsilon_{\mu}}\varepsilon_{\lambda\mu}'
z_{\lambda\mu}^b)^2).$

\noindent
For simplicity we write $\sqrt{\varepsilon_{\mu}}\varepsilon_{\lambda\mu}
'z_{\lambda\mu}^a$,
$\sqrt{\varepsilon_{\mu}}\varepsilon_{\lambda\mu}'z_{\lambda\mu}^b$ as
$a_{\lambda\mu}$, $b_{\lambda\mu}$. It is

$\exp(-\pi (x+ib_{\lambda\mu})^2)-\exp(-\pi
(a_{\lambda\mu}+ib_{\lambda\mu})^2)$

$=\exp(-\pi
(x^2-b_{\lambda\mu}^2))\exp(-2i\pi b_{\lambda\mu}x)-\exp(-\pi
(a_{\lambda\mu}^2-b_{\lambda\mu}^2))\exp(-2i\pi
b_{\lambda\mu}a_{\lambda\mu})$

$=\{\exp(-\pi (x^2-b_{\lambda\mu}^2))\cos(2\pi
b_{\lambda\mu}x)-\exp(-\pi (a_{\lambda\mu}^2-b_{\lambda\mu}^2)\}
\cos(2\pi b_{\lambda\mu}a_{\lambda\mu})\}$

$-i\{\exp(-\pi (x^2-b_{\lambda\mu}^2))\sin (2\pi b_{\lambda\mu}x)-\exp(-\pi
(a_{\lambda\mu}^2-b_{\lambda\mu}^2))
\sin (2\pi b_{\lambda\mu}a_{\lambda\mu})\}.\;\;\;\;\;\cdots(\ast_3)$

\noindent
We consider the first term of $(\ast_3)$. Then

$\exp(-\pi (x^2-b_{\lambda\mu}^2))\cos(2\pi b_{\lambda\mu}x)=\exp(\pi
b_{\lambda\mu}^2)\exp(-\pi x^2)\cos(2\pi b_{\lambda\mu}x).$

\noindent
We put $f(x)=\exp(-\pi x^2)\cos(2\pi b_{\lambda\mu}x).$ We assume that $0\le
b_{\lambda\mu}$.

$f'(x)=-2\pi x\exp(-\pi x^2)\cos(2\pi b_{\lambda\mu}x)-\exp(-\pi x^2)2\pi
b_{\lambda\mu}\sin(2\pi b_{\lambda\mu}x)$

$=-2\pi\sqrt{x^2+b_{\lambda\mu}^2}\exp(-\pi x^2)\cos(2\pi
b_{\lambda\mu}x+\alpha_x),$

\noindent
where $\cos\alpha_x=\frac{x}{\sqrt{x^2+b_{\lambda\mu}^2}},\;\; -\sin
\alpha_x=\frac{b_{\lambda\mu}}{\sqrt{x^2+b_{\lambda\mu}^2}}.$ Since
$-\sin \alpha_x=\frac{b_{\lambda\mu}}
{\sqrt{x^2+b_{\lambda\mu}^2}}$, $\alpha_x$ is negative. There is a unique
maximum of $|f(x)|$ in 

\noindent
$\left\{x\in{\bf R}\,\left|\,\frac{\pi}{2}(2m-1)\le 2\pi
b_{\lambda\mu}x<\frac{\pi}{2}(2m+1)\right\}\right.$ for each $m\in{\bf Z}$, that
is, $x$ satisfies

$f'(x)=0,\; \frac{\pi}{2}(2m-1)\le 2\pi
b_{\lambda\mu}x<\frac{\pi}{2}(2m+1)\Longleftrightarrow 2\pi b_{\lambda\mu}
x+\alpha_x=\frac{\pi}{2}(2m-1).
\;\;\;\;\;\cdots (\ast_4)$

\noindent
We write the value of $x$ having the maximum of $|f(x)|$ in the interval as
$x_m$. On the other hand, we denote the value $\alpha_x$ at $x=A_{2m}$  by
$\alpha_{A_{2m}}$. Then

$A_{2m}=\frac{m-\frac{1}{2}-\frac{\alpha_{A_{2m}}}{\pi}}{2b_{\lambda\mu}}.$
 The maximum of $f(x)$ is $f(A_{2m})=\exp\left(-\pi x_m^2\right)\cos \left(m\pi
-\frac{\pi}{2} -\alpha_{A_{2m}} \right)$

\noindent
Hence $|f(A_{2m})|\le \exp\left(-\pi
\left(\frac{m-\frac{1}{2}-\frac{\alpha_{A_{2m}}}{\pi}}{2b_{\lambda\mu}}
\right)^2\right)
\le \exp\left(-\pi \left(\frac{m-\frac{1}{2}}{2b_{\lambda\mu}}\right)
^2\right).$ Then

$\Biggl|\sum_{\varepsilon_{\lambda\mu}'z_{\lambda\mu}^a \in
L_{\lambda\mu}'}\sqrt{\varepsilon_{\mu}}\varepsilon'{}_{\lambda\mu} \exp(-\pi
(\sqrt{\varepsilon_{\mu}}\varepsilon_{\lambda\mu}'z_{\lambda\mu}^a)^2-(\varepsilon_{\lambda\mu}'z_{\lambda\mu}^b)^2))
\cos(2\pi
\varepsilon_{\lambda\mu}'z_{\lambda\mu}^b\sqrt{\varepsilon_{\mu}}
\varepsilon_{\lambda\mu}'z_{\lambda\mu}^a)$

$-\int_{-\sqrt{\varepsilon_{\mu}}\frac{H_{\lambda\mu}'}{2}}^{\sqrt{\varepsilon_{\mu}}\frac{H_{\lambda\mu}'}{2}}\exp(-\pi
(x^2-b_{\lambda\mu}^2))\cos(2\pi b_{\lambda\mu}x)dx\Biggr|$

$=\exp(\pi
b_{\lambda\mu}^2)\Biggl|\sum_{\varepsilon_{\lambda\mu}'z_{\lambda\mu}^a \in
L_{\lambda\mu}'}\sqrt{\varepsilon_{\mu}}\varepsilon'{}_{\lambda\mu}
f(\sqrt{\varepsilon_{\mu}}\varepsilon_{\lambda\mu}'z_{\lambda\mu}^a)^2)
-\int_{-\sqrt{\varepsilon_{\mu}}\frac{H_{\lambda\mu}'}{2}}^{\sqrt{
\varepsilon_{\mu}}\frac{H_{\lambda\mu}'}{2}}f(x)dx\Biggr|$.

\noindent
We denote $A_{2m+1}=\frac{1}{4b_{\lambda\mu}}(2m+1).$
Since $[(b(k)_{\lambda\mu})]$ is finite, there exists a positive real number 
$c$ so that $\star(\ast(c))\le \left[\left(\left|\frac{1}
{4b_{\lambda\mu}}\right|\right)\right],$
that is, $\left\{\lambda\,\left|\, \left\{\mu\,\left|\,c\le\left|
\frac{1}{4b_{\lambda\mu}}\right|\right\}\in F_2\right\}\in
F_1\right.\right..$
Furthermore since $\sqrt{\varepsilon}\varepsilon'$ is
infinitesimal in $\,^{\star}(\,^{\ast}{\bf R})$,
$k\in L \Rightarrow \sqrt{\varepsilon}\varepsilon' <\left|
\frac{1}{4b(k)}\right|,$
that is,
$\left\{\lambda\,\left|\, \left\{\mu\,\left|\sqrt{\varepsilon_{\mu}}
\varepsilon_{\lambda\mu}'
<\left|\frac{1}{4b_{\lambda\mu}}(k_{\mu})\right|\right\}\in
F_2\right\}\in F_1\right.\right..$
We assume $\lambda\mu$ satisfies the above condition.
We denote $\sqrt{\varepsilon_{\mu}}\varepsilon
_{\lambda\mu}'$ by $\Delta '$. We shall show that the following term is
infinitesimal in $\,^{\star}(\,^{\ast}{\bf R})$ :
$\bigl|\sum_{\varepsilon_{\lambda\mu}'z_{\lambda\mu}^a \in
L_{\lambda\mu}'}\sqrt{\varepsilon_{\mu}}\varepsilon'{}_{\lambda\mu}
f(\sqrt{\varepsilon_{\mu}}\varepsilon_{\lambda\mu}'z_{\lambda\mu}^a)^2)
-\int_{-\sqrt{\varepsilon_{\mu}}\frac{H_{\lambda\mu}'}{2}}^{\sqrt
{\varepsilon_{\mu}}\frac{H_{\lambda\mu}'}{2}}f(x)dx\bigr|.$
We write the maximum of $j$ as l so that $A_j\in [0,\sqrt
{\varepsilon_{\mu}}\frac{H_{\lambda\mu}'}{2}]$.
Let $x_k=\Delta ' k,\;\; x_{i_j}\le A_j<x_{i_j+1}\;(1\le j\le l)$. Then we 
devide the interval $[0,\sqrt
{\varepsilon_{\mu}}\frac{H_{\lambda\mu}'}{2}]$ into suitable intervals and prove
it.

\medskip

\textsc{Case} 1.

$\bigl|\int_{-A_1}^{A_1}f(x)dx-\left(\sum_{i=-i_1}^{i_1}f(x_i)\right)
\Delta ' \bigr|$

$=\bigl|\int_{-A_1}^{-x_{i_1}}f(x)dx+\sum_{i=-i_1}^{i_1}\left(
\int_{x_i}^{x_{i+1}}(f(x)-f(x_i))dx\right)$

$+\int_{x_{i_1}}^{A_1}(f(x)-f(x_{i_1}))dx-f(x_{i_1})(x_{{i_1}+1}-A_1)\bigr|$

$=\bigl|\int_{-A_1}^{-x_{i_1}}f(x)dx+
\int_{-x_{i_1}}^{-x_{{i_1}-1}}(f(x)-f(-x_{{i_1}}))dx+\cdots
+\int_{-x_{1}}^{x_{0}}(f(x)-f(-x_{1}))dx$

$+\int_{x_0}^{x_{1}}(f(x)-f(x_0))dx+\cdots
+\int_{x_{i_1}-1}^{x_{i_1}}(f(x)-f(x_{{i_1}-1}))dx$

$+\int_{x_{i_1}}^{A_1}(f(x)-f(x_{i_1}))dx-f(x_{i_1})(x_{{i_1}+1}-A_1)\bigr|\textrm{,
where }x_0=0,$

$\le\int_{-A_1}^{-x_{i_1}}f(x)dx+
\int_{-x_{i_1}}^{-x_{{i_1}-1}}(f(x)-f(-x_{{i_1}}))dx+\cdots
+\int_{-x_{1}}^{x_{0}}(f(x)-f(-x_{1}))dx$

$+\int_{x_0}^{x_{1}}(f(x_0)-f(x))dx+\cdots
+\int_{x_{i_1}-1}^{x_{i_1}}(f(x_{{i_1}-1})-f(x))dx$

$+\int_{x_{i_1}}^{A_1}(f(x_{i_1})-f(x))dx+f(x_{i_1})(x_{{i_1}+1}-A_1)$

$\textrm{, since } f\textrm{ is an even function, }f(-x)=f(x),$

$=\int_{x_{i_1}}^{A_1}f(x)dx+
\int_{x_{{i_1}-1}}^{x_{i_1}}(f(x)-f(x_{{i_1}}))dx+\cdots
+\int_{-x_{0}}^{x_{1}}(f(x)-f(x_{1}))dx$

$+\int_{x_0}^{x_{1}}(f(x_0)-f(x))dx+\cdots
+\int_{x_{i_1}-1}^{x_{i_1}}(f(x_{{i_1}-1})-f(x))dx$

$+\int_{x_{i_1}}^{A_1}(f(x_{i_1})-f(x))dx+f(x_{i_1})(x_{{i_1}+1}-A_1)$

$=\left(f(x_0)-f(x_1)+f(x_1)-f(x_2)+\cdots
+f(x_{{i_1}-1})-f(x_{i_1})+f(x_{i_1})\right)\Delta '$

$=f(x_0)\Delta ' =f(0)\Delta '$.

\bigskip

\textsc{Case} 2.

$\bigl|\int_{-x_{i_2}}^{-A_1}f(x)dx+\int_{A_1}^{x_{{i_2}+1}}f(x)dx
-\left(\sum_{i=-{i_2}}^{-({i_1}+1)}f(x_i)+\sum_{i={i_1}+1}^{i_2}f(x_i)
\right)\Delta '\bigr|$

$=\bigl|\int_{-x_{{i_2}}}^{-x_{{i_2}-1}}(f(x)-f(-x_{{i_2}}))dx+\cdots +
\int_{-x_{{i_1}+2}}^{-x_{{i_1}+1}}(f(x)-f(-x_{{i_1}+2}))dx$

$+\int_{-x_{{i_1}+1}}^{-A_{1}}(f(x)-f(-x_{{i_1}+1}))dx
-f(-x_{{i_1}+1})(A_1-x_{i_1})+\int_{A_1}^{x_{{i_1}+1}}f(x)dx$

$+\int_{x_{{i_1}+1}}^{x_{{i_1}+2}}(f(x)-f(x_{{i_1}+1}))dx+\cdots +
\int_{x_{i_2}}^{x_{{i_2}+1}}(f(x)-f(x_{i_2}))dx\bigr|$

$\le\int_{-x_{{i_2}}}^{-x_{{i_2}-1}}(f(x)-f(-x_{{i_2}}))dx+\cdots +
\int_{-x_{{i_1}+2}}^{-x_{{i_1}+1}}(f(x)-f(-x_{{i_1}+2}))dx$

$+\int_{-x_{{i_1}+1}}^{-A_{1}}(f(x)-f(-x_{{i_1}+1}))dx
-f(-x_{{i_1}+1})(A_1-x_{i_1})+\int_{A_1}^{x_{{i_1}+1}}(-f(x))dx$

$+\int_{x_{{i_1}+1}}^{x_{{i_1}+2}}(f(x_{{i_1}+1})-f(x))dx+\cdots +
\int_{x_{i_2}}^{x_{{i_2}+1}}(f(x_{i_2})-f(x))dx$

$\textrm{, since }f\textrm{ is an even function, }f(-x)=f(x),$

$=\int_{x_{{i_2}-1}}^{x_{i_2}}(f(x)-f(x_{{i_2}}))dx+\cdots +
\int_{x_{{i_1}+1}}^{x_{{i_1}+2}}(f(x)-f(x_{{i_1}+2}))dx$

$+\int_{A_1}^{x_{{i_1}+1}}(f(x)-f(x_{{i_1}+1}))dx
-f(-x_{{i_1}+1})(A_1-x_{i_1})+\int_{A_1}^{x_{{i_1}+1}}(-f(x))dx$

$+\int_{x_{{i_1}+1}}^{x_{{i_1}+2}}(f(x_{{i_1}+1})-f(x))dx+\cdots +
\int_{x_{i_2}}^{x_{{i_2}+1}}(f(x_{i_2})-f(x))dx,$

$\textrm{,since }\int_{x_{i_2}}^{x_{{i_2}+1}}(f(x_{i_2})-f(x))dx\le
\int_{x_{i_2}}^{x_{{i_2}+1}}(f(x_{i_2})-f(A_2))dx,$

$\le(-f(x_{{i_1}+1})+f(x_{{i_1}+1})-f(x_{{i_1}+2})+\cdots
+f(x_{{i_2}-1})-f(x_{i_2})+f(x_{i_2})-f(A_2))\Delta '$

$=(-f(A_2))\Delta '.$

\medskip

Since the next steps are just same, we obtain the following :

$\bigl|\sum_{\varepsilon_{\lambda\mu}'z_{\lambda\mu}^a \in
L_{\lambda\mu}'}\sqrt{\varepsilon_{\mu}}\varepsilon'{}_{\lambda\mu}
f(\sqrt{\varepsilon_{\mu}}\varepsilon_{\lambda\mu}'z_{\lambda\mu}^a)^2)
-\int_{-\sqrt{\varepsilon_{\mu}}\frac{H_{\lambda\mu}'}{2}}^{\sqrt{
\varepsilon_{\mu}}\frac{H_{\lambda\mu}'}{2}}f(x)dx\bigr|$

$\le
2\sqrt{\varepsilon_{\mu}}\varepsilon'{}_{\lambda\mu}\sum_{m=0}^{l}\exp\left(-\pi
\left(\frac{m-\frac 12}{2b_{\lambda\mu}}\right)^2\right).$

\noindent
Hence if $0\le b_{\lambda\mu}$,

$\bigl|\sum_{\varepsilon_{\lambda\mu}'z_{\lambda\mu}^a,
\varepsilon_{\lambda\mu}'z_{\lambda\mu}^b \in
L_{\lambda\mu}'}\sqrt{\varepsilon_{\mu}}\varepsilon'{}_{\lambda\mu} \exp(-\pi
(\sqrt{\varepsilon_{\mu}}\varepsilon_{\lambda\mu}'z_{\lambda\mu}^a)^2-
(\varepsilon_{\lambda\mu}'z_{\lambda\mu}^b)^2))
\cos(2\pi
\varepsilon_{\lambda\mu}'z_{\lambda\mu}^b\sqrt{\varepsilon_{\mu}}
\varepsilon_{\lambda\mu}'z_{\lambda\mu}^a)$

$-\int_{-\sqrt{\varepsilon_{\mu}}\frac{H_{\lambda\mu}'}{2}}^{\sqrt{\varepsilon_{\mu}}\frac{H_{\lambda\mu}'}{2}}\exp(-\pi
(x^2-b_{\lambda\mu}^2))\cos(2\pi b_{\lambda\mu}x)dx\bigr|$

$\le
\exp(\pi
b_{\lambda\mu}^2)\cdot
2\sqrt{\varepsilon_{\mu}}\varepsilon'{}_{\lambda\mu}
\sum_{m=0}^{l}\exp\left(-\pi
\left(\frac{m-\frac 12}{2b_{\lambda\mu}}\right)^2\right).$

\noindent
Since
$\left[\left(\sum_{m=0}^{l}\exp\left(-\pi \left(\frac{m-\frac
12}{2b_{\lambda\mu}}\right)^2\right)\right)\right]\le 
\sum_{m=0}^{l}\exp\left(-\pi \left(\frac{m-\frac
12}{2b_0}\right)^2\right),$
it is finite. Hence
$\left[\left(\exp(\pi
b_{\lambda\mu}^2)\cdot
2\sqrt{\varepsilon_{\mu}}\varepsilon'{}_{\lambda\mu}\sum_{m=0}^{l}\exp\left(-\pi
\left(\frac{m-\frac 12}{2b_{\lambda\mu}}\right)^2\right)\right)\right]$
is infinitesimal in $\,^{\star}(\,^{\ast}{\bf R})$. If $b_{\lambda\mu}<0$, 
the argument is parallel, and also, for the term of sin in $(\ast_3)$, though
sin is not an even function, the same argument holds. Hence
$[(B_{\lambda\mu}(k_\mu))]$ is infinitesimal in $\,^{\star}(\,^{\ast}{\bf C})$
with respect to ${\bf C}$. Put $N_{\lambda\mu}=\frac{1}{B_{\lambda\mu}(k_\mu)}$.
Then it is  infinite, and since
$\lim_{n\to\infty}\left(1+\frac{1}{an}\right)^n=e^{-a},$

$\displaystyle\textrm{{\bf st}}\left(\left[\left(
1+\frac{1}{\int_{-\infty}^{\infty}\exp(-\pi x^2)dx\cdot \frac{1}{B_{\lambda\mu}(k_\mu)}}\right)^{\frac{1}{B_{\lambda\mu}(k_\mu)}}
\right]\right)$ 

$=\displaystyle\textrm{{\bf st}}\left(\left[\left(
1+\frac{1}{\int_{-\infty}^{\infty}\exp(-\pi x^2)dx\cdot N_{\lambda\mu}}\right)^{N_{\lambda\mu}}
\right]\right)=\exp\left(-\int_{-\infty}^{\infty}\exp(-\pi x^2)dx\right).$

\noindent
There is an infinitesimal $C\; (=[(C_{\lambda\mu})])$ in 
$\,^{\star}(\,^{\ast}{\bf C})$ so that

$\displaystyle\left(1+\frac{1}{\int_{-\infty}^{\infty}\exp(-\pi x^2)dx\cdot
N_{\lambda\mu}}\right)^{N_{\lambda\mu}}=\exp\left(-\int_{-\infty}^{\infty}\exp(-\pi
x^2)dx\right)+C_{\lambda\mu}.$

\noindent
Then $[((\ast_1))]$ is equal to

$\displaystyle\left[\left(\left(1+\frac{1}{\int_{-\infty}^{\infty}\exp(-\pi
x^2)dx\cdot
N_{\lambda\mu}}\right)^{N_{\lambda\mu}}\right)^{[(B_{\lambda\mu}(k_\mu)
H_{\mu}^2)]}\right]$

$=\left(\star\left(\ast\left(\exp\left(-\int_{-\infty}^{\infty}\exp(-\pi
x^2)dx\right)\right)\right)+C \right)^{[(B_{\lambda\mu}(k_\mu)H_{\mu}^2)]}.$

\noindent
Since $[(B_{\lambda\mu}(k_\mu)H_{\mu}^2)]$ is infinitesimal,

$\textrm{{\bf
st}}\left(\left(\star\left(\ast\left(\exp\left(-\int_{-\infty}^{\infty}\exp(-\pi
x^2)dx\right)\right)\right)+C\right)^{[(B_{\lambda\mu}(k_\mu)H_{\mu}^2)]}
\right)=1.$

\noindent
Thus

$\displaystyle\textrm{st}\left(\textrm{st}\left(\frac{\sum_{a\in
X}\varepsilon_0\exp
\left(-\pi\,^{\star}\varepsilon\sum_{k\in
L}(a(k)+ib(k))^2\right)}{\star\left(\ast\left(
\int_{-\infty}^{\infty}\exp(-\pi
x^2)dx\right)^{H^2}\right)}\right)\right)=1.$

\noindent
Since $\int_{-\infty}^{\infty}\exp(-\pi
x^2)dx=1$, then st(st($C_2(b)$))$=1$.

\bigskip 

{\it Acknowledgement}.  We would like to thank Prof. T.
Kamae for a useful suggestion about Gauss sums.

\bigskip 

\begin{center}
{\bf References}
\end{center}

[An]  \; R.M. Anderson, A non-standard representation for Brownian motion
and It\^o integration, Israel J. Math. {\bf 25} (1976), 15-46.

[A-F-HK-L]  \; S. Albeverio, J.E. Fenstad, R. H$\phi$egh-Krohn,
T. Lindstro$\phi$m, Non-standard methods in stochastic analysis and
mathematical physics, Academic Press (1986).

[F]  \; D. Fujiwara, A construction of the fundamental solution
for the Schr\"odinger equation, J. D'Analyse Math. {\bf 35} (1979),
41-96.

[F-H] \; R.P. Feynman, A.R. Hibbs, Quantum mechanics and path integrals,
McGrow-Hill Inc. All rights (1965).

[G] \; E.I. Gordon, Nonstandard methods in commutative harmonic analysis,
Translations of mathematical monographs {\bf 164} American mathematical
society, 1997.

[H]\; T. Hida, Brownian motion, in Japanese, Iwanami shoten, 1975.

[Ic]	\; T. Ichinose, Path integral for the Dirac equation in two space-time
dimensions, Proc. Japan Acad. Ser. A. Math. Sci. {\bf 58} (1982), 290-293.

[It]	\; T. It\^o, Differential equations determinig a Markoff process
(original Japanese: Zenkoku Sizyo Sugaku Danwakai-si), Journ. Pan-Japan Math.
Coll. No.{\bf 1077}, 1942.

[I-T]	\; T. Ichinose, H. Tamura, Path integral approach to relativistic
quantum mechanics-Two-dimensional Dirac equation,
Suppl. Prog. Theor. Phys. {\bf 92} (1987), 144-175.

[Ka]\; T. Kamae, A simple proof of the Ergodic theorem using non-standard
analysis, Israel J. Math. {\bf 42} (1982), 284-290.

[Ki]	\; M. Kinoshita, Nonstandard representation of distribution I, Osaka J.
Math. {\bf 25} (1988), 805-824.

[Loe] \; P.A. Loeb, Conversion from nonstandard to standard measure spaces and
application in probability theory, Trans. Amer. Math. Soc. {\bf 211} (1975),
113-122.

[Loo1] \; K. Loo, Nonstandard Feynman path integral for harmonic oscillator,
J. Math. Phys., {\bf 40} (1999), 5511-5521.

[Loo2] \; K. Loo, A rigourous real-time Feynman path integral and
propagator, J. Phys., {\bf A33} (2000), 9215-9239.

[Lu] \; W.A. Luxemburg, A Nonstandard approach to Fourier analysis,
Contributions to Nonstandard Analysis, North-Holland, Amsterdam, pp.16-39, 1972.

[Na1] \; T. Nakamura, A nonstandard reprentation of Feynman's Path
integrals, J. Math. Phys., {\bf 32} (1991), 457-463.

[Na2] \; T. Nakamura, Path space measure for the 3+1-dimensional Dirac
equation in momentum space, J. Math. Phys., {\bf 41} (2000), 5209-5222.

[Ne] \; E. Nelson, Feynman integrals and the Schr\"odinger equation, J.
Math. Phys., {\bf 5} (1964), 332-343.

[N-O]  \; T. Nitta and T. Okada, Double infinitesimal Fourier
transformation for the space of functionals and reformulation of Feynman path
integral, Lecture Note Series in Mathematics, Osaka University Vol.{\bf 7} (2002),
255-298 in
Japanese.

[N-O-T]\; T. Nitta, T. Okada and A. Tzouvaras, Classification of
non-well-founded sets and an application, Math. Log. Quart. {\bf 49}
(2003), 187-200.

[R] \; R. Remmert, Theory of complex functions, Graduate Texts in Mathematics
{\bf 122}, Springer, Berlin-Heidelberg-New York, 1992. 

[S]\; M. Saito, Ultraproduct and non-standard analysis, in Japanese, Tokyo
tosho, 1976.

[T]\; G. Takeuti, Dirac space, Proc. Japan Acad. {\bf 38}
(1962), 414-418.

\newpage

{\small Takashi NITTA

Department of Education

Mie University 

Kamihama, Tsu, 514-8507, Japan

e-mail : nitta@edu.mie-u.ac.jp

\bigskip

Tomoko OKADA

Graduate school of Mathematics

Nagoya University

Chikusa-ku, Nagoya, 464-8602, Japan

e-mail : m98122c@math.nagoya-u.ac.jp}

\end{document}